# ITERATED SMOOTHED BOOTSTRAP CONFIDENCE INTERVALS FOR POPULATION QUANTILES

By Yvonne H. S. Ho and Stephen M. S. Lee[1]

*The University of Hong Kong*

This paper investigates the effects of smoothed bootstrap iterations on coverage probabilities of smoothed bootstrap and bootstrap-$t$ confidence intervals for population quantiles, and establishes the optimal kernel bandwidths at various stages of the smoothing procedures. The conventional smoothed bootstrap and bootstrap-$t$ methods have been known to yield one-sided coverage errors of orders $O(n^{-1/2})$ and $o(n^{-2/3})$, respectively, for intervals based on the sample quantile of a random sample of size $n$. We sharpen the latter result to $O(n^{-5/6})$ with proper choices of bandwidths at the bootstrapping and Studentization steps. We show further that calibration of the nominal coverage level by means of the iterated bootstrap succeeds in reducing the coverage error of the smoothed bootstrap percentile interval to the order $O(n^{-2/3})$ and that of the smoothed bootstrap-$t$ interval to $O(n^{-58/57})$, provided that bandwidths are selected of appropriate orders. Simulation results confirm our asymptotic findings, suggesting that the iterated smoothed bootstrap-$t$ method yields the most accurate coverage. On the other hand, the iterated smoothed bootstrap percentile method interval has the advantage of being shorter and more stable than the bootstrap-$t$ intervals.

**1. Introduction.** It is generally known that under Bhattacharya and Ghosh's (1978) smooth function model, the bootstrap percentile method confidence interval is subject to a one-sided coverage error of order $O(n^{-1/2})$, rendering it indistinguishable from the classical normal approximation method. Hall (1986) shows that Studentization can be employed to reduce this error to $O(n^{-1})$. The iterated bootstrap offers an alternative to error correction by calibrating the nominal coverage level iteratively; see Beran (1987). Hall and

Received April 2003; revised March 2004.

[1]Supported by a grant from the Research Grants Council of the Hong Kong Special Administrative Region, China (Project HKU 7131/00P).

Martin (1988) show that each such iteration reduces the one-sided coverage error by an order of $O(n^{-1/2})$ successively. On the other hand, smoothing the bootstrap, which amounts to drawing bootstrap samples from a kernel-smoothed empirical distribution, instead of sampling with replacement from the raw data, does not affect the coverage accuracy of bootstrap intervals in general. Polansky and Schucany (1997) propose smoothed bootstrap strategies to yield intervals of $O(n^{-1})$ coverage error. Their methods, however, require sophisticated tuning of the smoothing bandwidths, rendering the improvement less stable than that resulting from Studentization or the iterated bootstrap.

In contrast to problems in the context of smooth function models, conventional bootstrap confidence intervals for the $q$th population quantile, for a fixed $q \in (0, 1)$, have notably poor coverages; see Hall and Martin (1989). Here the percentile method gives coverage error of order $O(n^{-1/2})$, which cannot be improved upon by nominal coverage calibration using the iterated bootstrap. Indeed, more generally, this $O(n^{-1/2})$ coverage error is inherent in any confidence interval procedure based directly on order statistics as a consequence of their binomial-type discreteness. See, for example, De Angelis, Hall and Young (1993) for a more detailed account of the above phenomenon. On the other hand, either the smoothed bootstrap or Studentization extends considerably the domain from which we derive the confidence limits, and may therefore be able to make asymptotic improvement over the conventional bootstrap percentile method.

In the context of estimating the variance $\sigma_n^2$ of the sample $q$th quantile, Hall, DiCiccio and Romano (1989) show that sufficiently high-order smoothing of the bootstrap succeeds in reducing the relative error of the unsmoothed bootstrap estimate from $O(n^{-1/4})$ to $O(n^{-1/2+\epsilon})$ for arbitrarily small $\epsilon > 0$. Falk and Janas (1992) show that smoothing the bootstrap percentile method returns the same order, $O(n^{-1/2})$, of coverage error as in the unsmoothed case. Studentization of the sample quantile involves estimation of $\sigma_n^2$, which admits an expansion $n^{-1}q(1-q)f(F^{-1}(q))^{-2} + O(n^{-3/2})$ under proper regularity conditions, where $f = F'$ and $F$ denotes the distribution function underlying the given random sample. In practice $\sigma_n^2$ has to be estimated from the sample by, for example, bootstrapping or explicit estimation of the leading term in its expansion above. Hall and Martin (1991) show that confidence intervals based on normal approximation of the sample quantile Studentized by the unsmoothed bootstrap variance estimate yield coverage error of order $O(n^{-1/2})$. Falk and Janas (1992) show that a similar result holds when the variance is estimated by plugging in the kernel density estimate of $f$. However, application of the smoothed bootstrap to the Studentized sample quantile in the latter case, which we shall term the smoothed bootstrap-$t$ method, succeeds in improving the error order to $o(n^{-2/3})$, if second-order kernels are used in conjunction with appropriately selected bandwidths at both the Studentization and bootstrapping steps; see Janas ([1993](#)). This result will be sharpened in Section 2.4, where we prove that the smoothed bootstrap-$t$ method can indeed yield a coverage error of precise order $O(n^{-5/6})$.

We investigate further in this paper the asymptotic effects on coverage error of iterating the smoothed bootstrap or bootstrap-$t$ methods, with the objective of generating confidence intervals with improved coverage accuracy. De Angelis, Hall and Young ([1993](#)) remark in passing without proof that iterating the smoothed bootstrap percentile method might reduce coverage error. They mention neither the degree of improvement nor the choices of kernel bandwidths, and implementation of the iterated smoothed bootstrap remains impractical. Iterating the smoothed bootstrap-$t$ has not been explored in the literature. We formalize the theory of the iterated smoothed bootstrap and bootstrap-$t$ methods by stating explicitly the orders of their coverage errors, and derive therefrom the optimal sizes of kernel bandwidths at the two levels of smoothed bootstrapping and also at the Studentization step where applicable. Specifically, we show for the smoothed bootstrap method that calibration of the nominal coverage level by one extra level of smoothed bootstrapping can reduce the coverage error from $O(n^{-1/2})$ to $O(n^{-2/3})$, provided that the bandwidths are chosen of order ranging from $n^{-2/9}$ to $n^{-1/12}$ at the outer level of bootstrapping, and of order $n^{-1/3}$ at the inner level. For the smoothed bootstrap-$t$ method, such coverage calibration succeeds in reducing the coverage error from $O(n^{-5/6})$ to $O(n^{-58/57})$, provided that the bandwidths are chosen of order $n^{-2/19}$ at the outer level, of order $n^{-11/57}$ at the inner level and of order ranging from $n^{-11/19}$ to $n^{-1/2}$ at the Studentization step. The latter result signifies thus far the best coverage accuracy achievable by bootstrap confidence intervals proposed in the literature for the population quantile. It also outperforms Beran and Hall's ([1993](#)) interpolated confidence interval, which is based on a convex combination of sample quantiles and yields a coverage error of order $O(n^{-1})$. Chen and Hall's ([1993](#)) smoothed empirical likelihood interval, which is based on smoothing the standard empirical likelihood procedure, has two-sided coverage error of order $O(n^{-1})$. They show also that Bartlett correction reduces the error further to $O(n^{-2})$ and a simple approximation to the correction results in error slightly smaller than $O(n^{-1})$. However, the one-sided coverage error of the smoothed empirical likelihood interval, Bartlett-corrected or not, remains of order $O(n^{-1/2})$. For, as we can see from (6.10) and (6.14) of Chen and Hall ([1993](#)), the one-sided coverage expansion for the interval consists of a nonvanishing term of order $n^{-1/2}$. This term persists even after Bartlett correction, which affects only terms of order $O(n^{-1})$. Our iterated smoothed bootstrap and bootstrap-$t$ intervals thus compare favorably with the smoothed empirical likelihood approach in terms of one-sided coverage accuracy. Moreover, our approaches have the further advantage of being extendable by additional bootstrap iterations to yield successively more accurate coverages, subject only to availability of computer resources and the extent to which asymptotic implications can be realized in practice.

Simulation results suggest that our two iterated bootstrap methods yield coverages much more accurate than their noniterated counterparts. The improvement is more significant in the non-Studentized case. In general our methods have accuracies comparable to the interpolation or smoothed empirical likelihood methods, but improving at a faster rate as $n$ increases. Despite its impressive coverage accuracy, the iterated smoothed bootstrap-$t$ method suffers, as expected, from the problem of instability pertinent to variance estimation for a sample quantile, which often results in overly long and highly variable confidence intervals. On the other hand, the iterated smoothed bootstrap percentile approach, albeit slightly less accurate, is relatively much more stable than both the noniterated and iterated smoothed bootstrap-$t$ methods.

Success of the iterated bootstrap in the present context extends its scope of application beyond the traditional regular problem settings and beyond the conventional, unsmoothed bootstrap procedures, yielding asymptotic improvement of a problem-specific nature. This confirms the potential of the iterated bootstrap as a general strategy for improving upon the bootstrap not just in cases where the conventional bootstrap works satisfactorily, such as in smooth function model settings, but also in cases where it does not work as satisfactorily, such as in the quantile case, and where a modified form of the bootstrap, such as the smoothed bootstrap, is required.

Section 2.1 introduces notation and states the regularity conditions required for the asymptotic theory. Sections 2.2 and 2.3 establish asymptotic expansions for the coverage probabilities of the noniterated and iterated smoothed bootstrap percentile method intervals, respectively, while their Studentized counterparts are treated in Sections 2.4 and 2.5. Based on our asymptotic results, we derive for each type of interval the optimal orders of kernel bandwidths at each level of bootstrapping and, where applicable, at the Studentization step in order to minimize coverage error. Section 2.6 discusses an alternative approach, which bases the confidence set root on a smoothed version of the sample quantile, to constructing bootstrap confidence intervals. Section 3 addresses the issue of empirical determination of bandwidths and suggests a bootstrap solution to the problem. Section 4 reports two simulation studies. The first demonstrates the bootstrap procedure for setting optimal bandwidths. The second shows that the iterated bootstrap improves upon the smoothed bootstrap and bootstrap-$t$ methods in terms of coverage accuracy. The iterated smoothed bootstrap interval also excels in terms of other indicators of interval performance. Section 5 concludes our findings. All technical details are given in the Appendix.

## 2. Theory.

2.1. *Notation and assumptions.* Let $\mathcal{X} = \{X_1, \ldots, X_n\}$ be independent and identically distributed from an unknown distribution $F$ with density $f = F'$. The parameter of interest is the $q$th population quantile $F^{-1}(q) \equiv \inf\{x \in \mathbb{R} : F(x) \geq q\}$, for a fixed $q \in (0, 1)$. We wish to construct a nominal level $1 - \alpha$ upper confidence interval for $F^{-1}(q)$, where $0 < \alpha < 1$.

Let $F_n$ be the empirical distribution of $\mathcal{X}$, which assigns a probability mass of $n^{-1}$ to each $X_i$, $i = 1, \ldots, n$, and let $\hat{F}_{n,\eta}$ be its kernel-smoothed version with density $\hat{f}_{n,\eta}$, such that

$$\hat{F}_{n,\eta}(t) = n^{-1} \sum_{i=1}^{n} K((t - X_i)/\eta) \quad \text{and} \quad \hat{f}_{n,\eta}(t) = (n\eta)^{-1} \sum_{i=1}^{n} k((t - X_i)/\eta),$$

for a kernel function $k$, $K(t) = \int_{-\infty}^{t} k(x)\,dx$ and a bandwidth $\eta > 0$. See Silverman ([1986](#)) for a general exposition of kernel density estimation. Note that the unsmoothed and smoothed sample $q$th quantiles are given, respectively, by $F_n^{-1}(q)$ and $\hat{F}_{n,\eta}^{-1}(q)$.

Write $\bar{q} = \min(q, 1 - q)$. Let $\Phi$ be the standard normal distribution function and let $\phi = \Phi'$ be its density. Define, for $\varepsilon \in (0, \bar{q})$, $D_1, D_2 > 0$ and $j = 1, 2, \ldots,$ $\mathcal{F}_j(\varepsilon, D_1, D_2)$ to be the class of distribution functions $\bar{F}$ satisfying the following smoothness conditions: (i) $\bar{F}^{-1}$ is $j$ times continuously differentiable in $(q - \varepsilon, q + \varepsilon)$, (ii) $(\bar{F}^{-1})'(q) \geq D_1$ and (iii) $\max_{i=1,\ldots,j} \sup_{p \in (q - \varepsilon, q + \varepsilon)} |(\bar{F}^{-1})^{(i)}(p)| \leq D_2$. It is clear that $\mathcal{F}_{j+1}(\varepsilon, D_1, D_2) \subset \mathcal{F}_j(\varepsilon, D_1, D_2)$ for $j = 1, 2, \ldots.$ We shall establish coverage expansions for our iterated and noniterated bootstrap intervals under $F \in \mathcal{F}_2(\varepsilon, D_1, D_2)$ for the non-Studentized case and $F \in \mathcal{F}_4(\varepsilon, D_1, D_2)$ for the Studentized case.

We make the following assumptions on the kernel $k$ throughout the paper:

(A1) $k$ is nonnegative, symmetric about zero and has compact support $[-a, a]$, for some $a > 0$;

(A2) $k^{(j)}$ exists and is bounded on $[-a, a]$ for $j = 1, 2, 3, 4$;

(A3) $\int_{-\infty}^{\infty} k(x)\,dx = \int_{-\infty}^{\infty} x^2 k(x)\,dx = 1$.

Note that the above assumptions require that $k$ be a proper density function, symmetric about 0, on the interval $[-a, a]$.

2.2. *Smoothed bootstrap percentile method.* Let $\mathcal{X}^{\dagger} = \{X_1^{\dagger}, \ldots, X_n^{\dagger}\}$ be a random sample simulated from $\hat{F}_{n,\eta}$. We may set in practice $X_i^{\dagger} = Y_i^* + \eta W_i^*$, $i = 1, \ldots, n$, where the $Y_i^*$ and the $W_i^*$ are independent random numbers drawn from the distributions $F_n$ and $K$, respectively. Denote by $F_{n,\eta}^*$ the empirical distribution of $\mathcal{X}^{\dagger}$. Define

$$G_n(t) = \mathbb{P}\{n^{1/2}(F_n^{-1}(q) - F^{-1}(q)) \leq t\}, \qquad t \in \mathbb{R}.$$

The smoothed bootstrap version of $G_n(t)$ substitutes $\hat{F}_{n,\eta}$ for $F$, and is given by

$$\hat{G}_{n,\eta}(t) = \mathbb{P}\{n^{1/2}(F_{n,\eta}^{*-1}(q) - \hat{F}_{n,\eta}^{-1}(q)) \le t | \mathcal{X}\}.$$

Following Beran's (1987) prepivoting idea, the root $n^{1/2}(F_n^{-1}(q) - F^{-1}(q))$ can be transformed, by prepivoting with $\hat{G}_{n,\eta}$, into a random variable approximately uniformly distributed over $[0, 1]$. It is clear in the context of confidence interval construction that the above action of prepivoting amounts to smoothed bootstrap estimation of the quantile of $n^{1/2}(F_n^{-1}(q) - F^{-1}(q))$. This defines a noniterated smoothed bootstrap percentile upper confidence interval, of nominal level $1 - \alpha$, to be

$$I_{1,\alpha} = (-\infty, F_n^{-1}(q) - n^{-1/2}\hat{G}_{n,\eta}^{-1}(\alpha)].$$

Write for brevity $\sigma_q = \{q(1-q)\}^{1/2}$. The following theorem establishes an asymptotic expansion for the distribution of the prepivoted root, and hence the coverage probability of $I_{1,\alpha}$.

THEOREM 1. *Under conditions* (A1)–(A3) *and assuming that* $\eta \propto n^{-\Delta_\eta}$ *with* $1/4 < \Delta_\eta < 1/2$, *we have that, for* $\alpha \in (0, 1)$,

$$
\begin{aligned}
(1) \quad & \mathbb{P}\{\hat{G}_{n,\eta}(n^{1/2}[F_n^{-1}(q) - F^{-1}(q)]) \ge \alpha\} \\
& \quad = \mathbb{P}\{F^{-1}(q) \in I_{1,\alpha}\} \\
& \quad = 1 - \alpha + n^{-1/2}\left(\frac{2q-1}{2\sigma_q} + \frac{\sigma_q f'(F^{-1}(q))}{f^2(F^{-1}(q))}\right)\Phi^{-1}(\alpha)^2 \phi(\Phi^{-1}(\alpha)) \\
& \qquad + O(\eta^2 + n^{-2}\eta^{-4} + n^{-1/2}\eta^{1/2})
\end{aligned}
$$

*uniformly in* $F \in \mathcal{F}_2(\varepsilon, D_1, D_2)$, *for any* $\varepsilon \in (0, \bar{q})$ *and* $D_1, D_2 > 0$.

We see from Theorem 1 that $I_{1,\alpha}$ has coverage error of precise order $O(n^{-1/2})$, provided that $F \in \mathcal{F}_2(\varepsilon, D_1, D_2)$ and the bandwidth $\eta \propto n^{-\Delta_\eta}$ is chosen such that $1/4 < \Delta_\eta \le 3/8$. Falk and Janas (1992) obtain an expansion similar to (1) for the coverage probability, up to order $o(n^{-1/2})$, under the restrictive condition that $\eta = o(n^{-1/3})$. The expansion (1) in our Theorem 1 gives an error term that reveals the explicit influence of the bandwidth $\eta$ on the coverage, which is crucial to our subsequent study of the effects of bootstrap iterations.

2.3. *Iterating the smoothed bootstrap percentile method.* In standard situations, the iterated bootstrap has been known to be very effective in enhancing coverage accuracy of confidence intervals. It operates by calibrating either the nominal level or the interval end points, with the use of an additional level of bootstrapping. We shall focus only on the former approach, which conforms exactly to Beran's (1987) prepivoting idea.

Define $\hat{G}_{n,\eta}$ as in Section 2.2 based on a generic smoothed bootstrap sample $\mathcal{X}^{\dagger}$ drawn from $\hat{F}_{n,\eta}$. Let $\mathcal{X}^{*} = \{X_1^{*}, \ldots, X_n^{*}\}$ be a generic outer-level smoothed bootstrap sample drawn from $\hat{F}_{n,\beta}$, for a kernel bandwidth $\beta > 0$, and let $F_{n,\beta}^{*}$ be its empirical distribution. Define, for $\eta > 0$, $\hat{H}_{n,\eta}$ by $\hat{H}_{n,\eta}(t) = n^{-1} \sum_{i=1}^{n} K((t - X_i^{*})/\eta)$, which corresponds to a smoothed empirical distribution of $\mathcal{X}^{*}$. Similarly we denote by $\mathcal{X}^{**}$ a generic inner-level smoothed bootstrap sample drawn from $\hat{H}_{n,\eta}$, and by $H_{n,\eta}^{*}$ the (unsmoothed) empirical distribution of $\mathcal{X}^{**}$. Define

$$\hat{G}_{n,\eta}^{*}(t) = \mathbb{P}\{n^{1/2}(H_{n,\eta}^{*-1}(q) - \hat{H}_{n,\eta}^{-1}(q)) \leq t | \mathcal{X}, \mathcal{X}^{*}\}, \qquad t \in \mathbb{R}.$$

Then the smoothed bootstrap estimates the distribution function of the prepivoted root $\hat{G}_{n,\eta}(n^{1/2}[F_n^{-1}(q) - F^{-1}(q)])$ by the conditional distribution, $\hat{J}_{n,\beta,\eta}$ say, of $\hat{G}_{n,\eta}^{*}(n^{1/2}[F_{n,\beta}^{*-1}(q) - \hat{F}_{n,\beta}^{-1}(q)])$ given $\mathcal{X}$. Prepivoting with $\hat{J}_{n,\beta,\eta}$ leads to a twice-prepivoted root $\hat{J}_{n,\beta,\eta}(\hat{G}_{n,\eta}(n^{1/2}[F_n^{-1}(q) - F^{-1}(q)]))$, which is asymptotically uniformly distributed over $[0, 1]$. In the context of confidence interval construction, this amounts to estimation of the $\alpha$th quantile of $n^{1/2}(F_n^{-1}(q) - F^{-1}(q))$ by $\hat{G}_{n,\eta}^{-1}(\hat{J}_{n,\beta,\eta}^{-1}(\alpha))$. The corresponding iterated smoothed bootstrap upper confidence interval, of nominal level $1 - \alpha$, is

$$I_{2,\alpha} = (-\infty, F_n^{-1}(q) - n^{-1/2} \hat{G}_{n,\eta}^{-1}(\hat{J}_{n,\beta,\eta}^{-1}(\alpha))].$$

Note that in the above procedure we have allowed use of two different bandwidths, $\beta$ and $\eta$, for the two levels of smoothed bootstrapping. This proves to be crucial to achieving asymptotic improvement in terms of coverage accuracy by means of the iterated bootstrap. The following theorem states in asymptotic terms how close the twice-prepivoted root is to a uniform random variable, and hence establishes the coverage error of $I_{2,\alpha}$.

THEOREM 2. *Assume the conditions in Theorem* 1, *that* $F \in \mathcal{F}_2(\varepsilon, D_1, D_2)$ *for some* $\varepsilon \in (0, \bar{q})$ *and* $D_1, D_2 > 0$, *and that* $\beta \propto n^{-\Delta_\beta}$ *with* $0 < \Delta_\beta < 1/3$. *Then we have*

$$\mathbb{P}\{\hat{J}_{n,\beta,\eta}(\hat{G}_{n,\eta}(n^{1/2}[F_n^{-1}(q) - F^{-1}(q)])) \geq \alpha\}$$
$$= \mathbb{P}\{F^{-1}(q) \in I_{2,\alpha}\}$$
$$= 1 - \alpha + O(\eta^2 + n^{-2}\eta^{-4} + n^{-1/2}\eta^{1/2} + n^{-1/2}\beta^2 + n^{-1}\beta^{-3/2}).$$

We see from Theorem 2 that the two levels of smoothed bootstrapping contribute separately to the coverage error of $I_{2,\alpha}$, which can be minimized to achieve $O(n^{-2/3})$ by setting $\eta \propto n^{-1/3}$ and $\beta \propto n^{-\Delta_\beta}$ with $1/12 \leq \Delta_\beta \leq 2/9$. The iterated smoothed bootstrap method thus improves upon the noniterated $I_{1,\alpha}$. We note, however, that application of the iterated smoothed bootstrap in the quantile problem does not yield the same level of improvement as has been well known in smooth function model situations, where each iteration of the (unsmoothed) bootstrap reduces coverage error by an order of $O(n^{-1/2})$.

2.4. *Smoothed bootstrap-t method.* We review in this section the smoothed bootstrap-$t$ method and derive explicitly the optimal orders of bandwidths that minimize its coverage error. Janas (1993) establishes that the smoothed bootstrap-$t$ method yields coverage error of order $o(n^{-2/3})$. Our results sharpen those of Janas (1993) by giving the precise order, namely $O(n^{-5/6})$, of the minimum coverage error. Noting that

$$\sigma_n^2 = \mathrm{Var}(F_n^{-1}(q)) = n^{-1}q(1-q)f(F^{-1}(q))^{-2} + O(n^{-3/2}),$$

we may estimate $n\sigma_n^2$, on plugging in a kernel density estimate for $f$, by

$$\hat{s}_\xi^2 = q(1-q)(n\xi)^2 \left\{ \sum_{i=1}^{n} h((F_n^{-1}(q) - X_i)/\xi) \right\}^{-2},$$

for some bandwidth $\xi > 0$ and kernel function $h$, which is assumed to satisfy:

(B1)   $h$ is nonnegative, symmetric about zero and has a compact support $[-b, b]$, for some $b > 0$;

(B2)   for some decomposition $-b = x_0 < x_1 < \cdots < x_m = b$, $h'$ exists on each interval $(x_{j-1}, x_j)$, is bounded and is either strictly positive or strictly negative there;

(B3)   $\int_{-\infty}^{\infty} h(x)\,dx = \int_{-\infty}^{\infty} x^2 h(x)\,dx = 1$.

Recall that $\mathcal{X}^\dagger = \{X_1^\dagger, \ldots, X_n^\dagger\}$ denotes a random sample from $\hat{F}_{n,\eta}$. Then the smoothed bootstrap version of $\hat{s}_\xi^2$ is given by

$$\hat{s}_\xi^{\dagger 2} = q(1-q)(n\xi)^2 \left\{ \sum_{i=1}^{n} h((F_{n,\eta}^{*-1}(q) - X_i^\dagger)/\xi) \right\}^{-2}.$$

Define, for $t \in \mathbb{R}$,

$$K_{n,\xi}(t) = \mathbb{P}\{n^{1/2}(F_n^{-1}(q) - F^{-1}(q))/\hat{s}_\xi \le t\}$$

and

$$\hat{K}_{n,\eta,\xi}(t) = \mathbb{P}\{n^{1/2}(F_{n,\eta}^{*-1}(q) - \hat{F}_{n,\eta}^{-1}(q))/\hat{s}_\xi^\dagger \le t | \mathcal{X}\}.$$

Then the smoothed bootstrap-$t$ upper confidence interval of nominal level $1 - \alpha$ is

$$I_{3,\alpha} = (-\infty, F_n^{-1}(q) - n^{-1/2}\hat{s}_\xi \hat{K}_{n,\eta,\xi}^{-1}(\alpha)].$$

Janas (1993) shows that if both $h$ and $k$ are chosen to be second-order, $\beta \propto n^{-1/3}$ and $\eta \propto n^{-1/5}$, the coverage error of $I_{3,\alpha}$ achieves an order of $o(n^{-2/3})$. We shall show that the minimum order of this coverage error is in fact $O(n^{-5/6})$, provided that the orders of bandwidths $\eta, \xi$ are chosen properly.

Similar to Theorem 1, the following theorem establishes an asymptotic expansion for the distribution of the prepivoted root $\hat{K}_{n,\eta,\xi}(n^{1/2}[(F_n^{-1}(q) - F^{-1}(q))/\hat{s}_\xi])$, and hence the coverage probability of $I_{3,\alpha}$.

THEOREM 3. *Assume conditions* (A1)–(A3), (B1)–(B3), *and that* $\eta \propto n^{-\Delta_\eta}$ *and* $\xi \propto n^{-\Delta_\xi}$, *with* $0 < \Delta_\eta < 1/5 < \Delta_\xi < 1$. *Then, for* $\alpha \in (0,1)$,

$$
\begin{aligned}
\mathbb{P}\{\hat{K}&_{n,\eta,\xi}(n^{1/2}[(F_n^{-1}(q) - F^{-1}(q))/\hat{s}_\xi]) \geq \alpha\} \\
&= \mathbb{P}\{F^{-1}(q) \in I_{3,\alpha}\} \\
&= 1 - \alpha + (n\eta)^{-1}D_{1,1}(F) + (n\xi)^{-1}\eta^2 D_{2,2}(F) \\
&\quad + n^{-1/2}\eta^2 D_{3,3}(F) + n^{-1}D_{4,2}(F) + n^{-3/2}\xi^{-1}D_{5,1}(F) \\
&\quad + O((n\xi)^{-5/2} + n\xi^5 + \xi^{5/2} + n^{-1/2}\eta^4 \\
&\qquad + n^{-1}\eta^{1/2} + n^{-3/2}\eta^{-5/2} + n^{-3/2}\xi^{-1/2}\eta^{-1} \\
&\qquad + (n\xi)^{-1}\eta^4 + n^{-1}\xi^{-1/2}\eta^2 + (n\eta)^{-1}\xi^{1/2} \\
&\qquad + n^{-1/2}\xi^{1/2}\eta^2 + n^{-3/2}\xi^{-1}\eta^{1/2} + n^{-1/2}\xi^2\eta^{-5/2} + (n\xi)^{-3/2}\eta^2)
\end{aligned}
$$

(2)

*uniformly in* $F \in \mathcal{F}_4(\varepsilon, D_1, D_2)$, *for any* $\varepsilon \in (0, \bar{q})$ *and* $D_1, D_2 > 0$. *Here, for each* $j = 1, \ldots, 5$, $D_{j,i}(F)$ *denotes a smooth function of the density derivatives* $\{f(F^{-1}(q)), f'(F^{-1}(q)), \ldots, f^{(i)}(F^{-1}(q))\}$.

We see from Theorem 3 that $I_{3,\alpha}$ has coverage error of order $O(n^{-5/6})$, provided that $\Delta_\eta = 1/6$ and $3/8 \leq \Delta_\xi \leq 1/2$. This suggests that $I_{3,\alpha}$ is asymptotically more advantageous than both $I_{1,\alpha}$ and $I_{2,\alpha}$. Note the different choices of bandwidth orders here as compared to Janas' (1993) recommendation, which yields only an $o(n^{-2/3})$ coverage error for $I_{3,\alpha}$.

It may be possible to further reduce the coverage error of either the iterated smoothed bootstrap or the smoothed bootstrap-$t$ method if a higher-order kernel function $k$ is employed. In fact, Janas (1993) shows that the error of $I_{3,\alpha}$ can be made as small as $O(n^{-1+\epsilon})$, for any $\epsilon > 0$, based on kernels of sufficiently high order. Similar results hold for the iterated smoothed bootstrap, suggesting that both $I_{2,\alpha}$ and $I_{3,\alpha}$ are essentially indistinguishable in terms of asymptotic coverage accuracy when high-order kernels are used. Our discussion nevertheless focuses on the practically more important second-order kernels, which have the virtue of being nonnegative and therefore allow straightforward Monte Carlo simulation from the resulting smoothed empirical distributions.

Studentization by $\hat{s}_\xi$ requires no direct simulation from the kernel $h$. We may therefore relax the second-order condition on $h$ without inflicting extra computational difficulty. Hall, DiCiccio and Romano (1989) show that use of a higher-order kernel $h$ can actually speed up the convergence rate of $\hat{s}_\xi$. However, we see from (2) that the best achievable coverage error of $I_{3,\alpha}$ is determined critically by the $(n\eta)^{-1}$ and $n^{-1/2}\eta^2$ terms. Increasing the order of $h$ affects only terms involving its bandwidth $\xi$, and can therefore not reduce the coverage error further.

It would be intriguing to explore the possibility of iterating the smoothed bootstrap-$t$ method to produce even more accurate confidence intervals for quantiles. We address this issue in the next section.

2.5. *Iterating the smoothed bootstrap-t method.* We follow the definitions used in Section 2.3 for bootstrap samples $\mathcal{X}^\dagger$, $\mathcal{X}^*$ and $\mathcal{X}^{**}$. The iterated smoothed bootstrap version of $\hat{s}^2$ is given by

$$\hat{s}_\xi^{**2} = q(1-q)(n\xi)^2 \left\{ \sum_{i=1}^n h((H_{n,\eta}^{*-1}(q) - X_i^{**})/\xi) \right\}^{-2}.$$

Define

$$\hat{K}_{n,\eta,\xi}^*(t) = \mathbb{P}\{n^{1/2}(H_{n,\eta}^{*-1}(q) - \hat{H}_{n,\eta}^{-1}(q))/\hat{s}_\xi^{**} \le t | \mathcal{X}, \mathcal{X}^*\}, \qquad t \in \mathbb{R}.$$

Similar to the construction of the iterated smoothed bootstrap interval $I_{2,\alpha}$, we set $\hat{L}_{n,\beta,\eta,\xi}$ to be the conditional distribution of $\hat{K}_{n,\eta,\xi}^*(n^{1/2}[F_{n,\beta}^{*-1}(q) - \hat{F}_{n,\beta}^{-1}(q)]/\hat{s}_\xi^*)$ given $\mathcal{X}$, where

$$\hat{s}_\xi^{*2} = q(1-q)(n\xi)^2 \left\{ \sum_{i=1}^n h((F_{n,\beta}^{*-1}(q) - X_i^*)/\xi) \right\}^{-2}.$$

The required twice-prepivoted root is given by

$$\hat{L}_{n,\beta,\eta,\xi}(\hat{K}_{n,\eta,\xi}(n^{1/2}[F_n^{-1}(q) - F^{-1}(q)]/\hat{s}_\xi)),$$

and the $\alpha$th quantile of $n^{1/2}(F_n^{-1}(q) - F^{-1}(q))/\hat{s}_\xi$ is estimated by $\hat{K}_{n,\eta,\xi}^{-1}(\hat{L}_{n,\beta,\eta,\xi}^{-1}(\alpha))$, where $\hat{K}_{n,\eta,\xi}$ and $\hat{s}_\xi$ are defined as in Section 2.4. The corresponding iterated smoothed bootstrap-$t$ upper confidence interval, of nominal level $1 - \alpha$, is

$$I_{4,\alpha} = (-\infty, F_n^{-1}(q) - n^{-1/2}\hat{s}_\xi \hat{K}_{n,\eta,\xi}^{-1}(\hat{L}_{n,\beta,\eta,\xi}^{-1}(\alpha))].$$

Note that construction of $I_{4,\alpha}$ involves three different bandwidths: $\xi$ at the Studentization step, $\beta$ at the outer level and $\eta$ at the inner level of smoothed bootstrapping. The following theorem establishes the order of the coverage error of $I_{4,\alpha}$ in terms of the three bandwidths.

TABLE 1
*Optimal bandwidth orders, on the $-\log_n$ scale, and corresponding coverage errors of $I_{j,\alpha}$, for $j = 1, 2, 3, 4$*

| | $-\log_n$(**bandwidth**) | | | |
| | $I_{1,\alpha}$ | $I_{2,\alpha}$ | $I_{3,\alpha}$ | $I_{4,\alpha}$ |
|---|---|---|---|---|
| Outer-level | — | $[1/12, 2/9]$ | — | $2/19$ |
| Inner-level | $(1/4, 3/8)$ | $1/3$ | $1/6$ | $11/57$ |
| Studentization step | — | — | $[3/8, 1/2)$ | $[1/2, 11/19)$ |
| Coverage error | $O(n^{-1/2})$ | $O(n^{-2/3})$ | $O(n^{-5/6})$ | $O(n^{-58/57})$ |

THEOREM 4. *Assume the conditions of Theorem 3, that $F \in \mathcal{F}_4(\varepsilon, D_1, D_2)$ for some $\varepsilon \in (0, \bar{q})$ and $D_1, D_2 > 0$, and that $\beta \propto n^{-\Delta_\beta}$ with $0 < \Delta_\beta < 1/7$. Then*

$$
\begin{aligned}
\mathbb{P}\{&\hat{L}_{n,\beta,\eta,\xi}(\hat{K}_{n,\eta,\xi}(n^{1/2}[F_n^{-1}(q) - F^{-1}(q)]/\hat{s}_\xi)) \geq \alpha\} \\
&= \mathbb{P}\{F^{-1}(q) \in I_{4,\alpha}\} \\
&= 1 - \alpha + O([(n\eta)^{-1} + n^{-3/2}\xi^{-1}](\beta^2 + n^{-1/2}\beta^{-3/2}) \\
&\qquad + [(n\xi)^{-1}\eta^2 + n^{-1}](\beta^2 + n^{-1/2}\beta^{-5/2}) \\
&\qquad + n^{-1/2}\eta^2(\beta^2 + n^{-1/2}\beta^{-7/2}) \\
&\qquad + (n\xi)^{-5/2} + n\xi^5 + \xi^{5/2} + n^{-1/2}\eta^4 \\
&\qquad + n^{-1}\eta^{1/2} + n^{-3/2}\eta^{-5/2} + n^{-3/2}\xi^{-1/2}\eta^{-1} \\
&\qquad + (n\xi)^{-1}\eta^4 + n^{-1}\xi^{-1/2}\eta^2 \\
&\qquad + (n\eta)^{-1}\xi^{1/2} + n^{-1/2}\xi^{1/2}\eta^2 + n^{-3/2}\xi^{-1}\eta^{1/2} \\
&\qquad\qquad + n^{-1/2}\xi^2\eta^{-5/2} + (n\xi)^{-3/2}\eta^2).
\end{aligned}
$$

(3)

We see from Theorem 4 that the iterated smoothed bootstrap-$t$ interval $I_{4,\alpha}$ can achieve an $o(n^{-1})$ coverage error. For instance, setting $\Delta_\eta = 11/57$, $1/2 \leq \Delta_\xi \leq 11/19$ and $\Delta_\beta = 2/19$ guarantees a coverage error of order $O(n^{-58/57})$. The precise minimum order of coverage error of $I_{4,\alpha}$ can be derived from a more detailed but uninspiring asymptotic expansion than that displayed in (3), which we omit here for simplicity.

We remark that the above iterated smoothed bootstrap-$t$ construction gives the fastest convergence rate of coverage as compared to other, smoothed or unsmoothed, bootstrap constructions thus far proposed in the literature. Not even the use of a high-order kernel $k$, which typically yields a coverage error of order $O(n^{-1+\epsilon})$ for an arbitrarily small $\epsilon > 0$ and a sufficiently high kernel order, is able to surpass this result. Successive iterations of the smoothed bootstrap procedure reduce the coverage error of $I_{4,\alpha}$ further. The forbidding task of managing a large number of bandwidths in a single asymptotic expansion prevents us from exploring this option further, although we recognize its potential in producing asymptotic improvement. The interpolation method proposed by Beran and Hall (1993) gives a confidence interval of $O(n^{-1})$ coverage error, which cannot be improved upon by higher-order interpolations.

We see from (3) that the $O(n^{-58/57})$ coverage error of $I_{4,\alpha}$ is determined strictly by terms of orders $n^{-1}\eta^{-1}\beta^2, n^{-1}\eta^2\beta^{-7/2}$ and $n^{-3/2}\eta^{-5/2}$. Similar to the case of smoothed bootstrap-$t$ construction, the coverage error of $I_{4,\alpha}$ cannot be further reduced by increasing the order of the kernel function $h$ used for Studentization, which affects only terms involving $\xi$.

Table 1 above summarizes the optimal choices of bandwidth orders and the corresponding one-sided coverage errors for the four intervals $I_{1,\alpha}$, $I_{2,\alpha}$, $I_{3,\alpha}$ and $I_{4,\alpha}$.

### 2.6. *Smoothing the sample quantile*: *an alternative.* Define

$$\tilde{f}_{n,\zeta}(t) = (n\zeta)^{-1}\sum_{i=1}^{n}\kappa((t-X_i)/\zeta) \quad \text{and} \quad \tilde{F}_{n,\zeta}(t) = \int_{-\infty}^{t}\tilde{f}_{n,\zeta}(x)\,dx,$$

for a kernel function $\kappa$ and a bandwidth $\zeta > 0$. A smoothed version of the sample quantile $F_n^{-1}(q)$ may be defined by $\tilde{F}_{n,\zeta}^{-1}(q)$, which we term the smoothed sample $q$th quantile. We now examine the effects on coverage error of basing the bootstrap confidence intervals on the root $n^{1/2}(\tilde{F}_{n,\zeta}^{-1}(q) - F^{-1}(q))$ instead of $n^{1/2}(F_n^{-1}(q) - F^{-1}(q))$. Define

$$G_{n,\zeta}(t) = \mathbb{P}\{n^{1/2}(\tilde{F}_{n,\zeta}^{-1}(q) - F^{-1}(q)) \le t\}, \qquad t \in \mathbb{R},$$

and $\hat{G}_{n,\eta,\zeta}$ to be the smoothed bootstrap version of $G_{n,\zeta}$ with $F$ substituted by $\hat{F}_{n,\eta}$. The corresponding smoothed bootstrap percentile upper confidence interval, of nominal level $1 - \alpha$, is

$$I_{1,\alpha}^{\kappa} = (-\infty, \tilde{F}_{n,\zeta}^{-1}(q) - n^{-1/2}\hat{G}_{n,\eta,\zeta}^{-1}(\alpha)].$$

The following theorem establishes an asymptotic expansion for the coverage probability of $I_{1,\alpha}^{\kappa}$.

THEOREM 5. *Suppose that* $F \in \mathcal{F}_2(\varepsilon, D_1, D_2)$ *for some* $\varepsilon \in (0, \bar{q})$ *and* $D_1, D_2 > 0$, *and that* $\kappa$ *is a second-order nonnegative kernel function. Under conditions* (A1)–(A3) *and assuming that* $\eta \propto n^{-\Delta_\eta}$ *with* $1/4 < \Delta_\eta < 1/3$ *and* $\zeta \propto n^{-\Delta_\zeta}$ *with* $\Delta_\zeta > 3/8$, *we have, for* $\alpha \in (0, 1)$, *that*

$$(4) \quad \mathbb{P}\{F^{-1}(q) \in I_{1,\alpha}^{\kappa}\} = 1 - \alpha + n^{-1/2}E^{\kappa} + o(n^{-1/2}) + O(n^{3/2}\zeta^4 + \eta^{-3/2}\zeta^2),$$

*for some nontrivial constant* $E^{\kappa}$ *independent of* $n, \zeta$.

We see from Theorem 5 that $I_{1,\alpha}^{\kappa}$ has coverage error of precise order $n^{-1/2}$ provided that $\Delta_\zeta \geq 1/2$. In any case, the order of the coverage error cannot be reduced further by adjusting the bandwidth $\zeta$. The best achievable coverage errors of both $I_{1,\alpha}$ and $I_{1,\alpha}^{\kappa}$ are of order $O(n^{-1/2})$, so that basing the smoothed bootstrap interval on the smoothed sample quantile does not yield any asymptotic improvement upon that derived from the sample quantile. We conjecture that arguments similar to those proving Theorem 5 can be employed to show that the smoothed sample quantile approach has no advantage either in the Studentized case.

**3. Empirical determination of bandwidths.** We now turn to the problem of empirical determination of the optimal bandwidths in practice. Many different practical strategies have been proposed for bandwidth selection in other problem settings, which often permit natural adaptation to our present context. Plausible approaches include, for example, cross-validation, an extra level of bootstrapping and plugging-in of sample quantities into asymptotic expansions of optimal bandwidths, among others.

Despite its computational intensity, the bootstrap approach is arguably the most straightforward method for setting optimal bandwidths. A smoothed bootstrap procedure for setting bandwidths is as follows. First, we generate an outermost level of smoothed bootstrap samples from $\hat{F}_{n,\gamma}$, for some bandwidth $\gamma > 0$. The collection of such samples, denoted generically by $\mathcal{X}^\circ$, forms the basis for our estimation of coverage probabilities and hence the determination of the best bandwidths. With reference to the optimal orders of bandwidths provided in Table 1, we set up a grid of pilot values of bandwidths for use in our procedure. For example, when considering $I_{2,\alpha}$, we may select $\beta_r$'s to be evenly spaced points within the range $[An^{-2/9}, Bn^{-1/12}]$ and $\eta_s$'s to be evenly spaced points within the range $[Cn^{-1/3}, Dn^{-1/3}]$, for some $A, B > 0$ and $D > C > 0$. The outermost bandwidth $\gamma$ can be fixed to be some multiple, $M$ say with $M > 1$, of the largest pilot bandwidth attempted at the outer level of the smoothed bootstrap. In our example we can set $\gamma = MBn^{-1/12}$. This is in line with our perception that the parent sample is drawn from an underlying distribution smoother than the smoothed empirical distribution used for bootstrapping. For each combination $(\beta, \eta) = (\beta_r, \eta_s)$ and each sample $\mathcal{X}^\circ$, we construct $I_{2,\alpha}$ and check if $\hat{F}_{n,\gamma}^{-1}(q) \in I_{2,\alpha}$. The coverage probability of $I_{2,\alpha}$, for each bandwidth pair $(\beta_r, \eta_s)$, is estimated by averaging over all samples $\mathcal{X}^\circ$. The required bootstrap confidence interval $I_{2,\alpha}$ is then constructed using the pair of bandwidths $(\beta_r, \eta_s)$ that gives the minimum coverage error as estimated above. We note that this procedure explicitly ensures that the bandwidths selected have the optimal asymptotic orders as displayed in Table 1. Selection of bandwidths for construction of the other three bootstrap intervals can be dealt with in a similar way.

Table 2
*Estimated coverage probabilities of $I_{j,\alpha}$, $j = 1, 3$, with bandwidths determined by smoothed bootstrap, for $q = 0.5$ and $\alpha = 0.05, 0.10, 0.90, 0.95$*

| Interval | $1 - \alpha = 0.05$ | $1 - \alpha = 0.10$ | $1 - \alpha = 0.90$ | $1 - \alpha = 0.95$ |
|---|---|---|---|---|
| | Standard normal data, $N(0, 1)$ | | | |
| $I_{1,\alpha}$ | 0.057 | 0.137 | 0.957 | 0.955 |
| $I_{3,\alpha}$ | 0.051 | 0.079 | 0.943 | 0.968 |
| | Double-exponential data, $\frac{1}{2} \exp(-|x|)$ | | | |
| $I_{1,\alpha}$ | 0.035 | 0.101 | 0.974 | 0.948 |
| $I_{3,\alpha}$ | 0.058 | 0.112 | 0.904 | 0.955 |
| | Lognormal data, $\exp\{N(0, 1)\}$ | | | |
| $I_{1,\alpha}$ | 0.053 | 0.104 | 0.966 | 0.972 |
| $I_{3,\alpha}$ | 0.071 | 0.094 | 0.809 | 0.893 |

**4. Simulation studies.** Two simulation studies were conducted to investigate the finite-sample performances of our proposed intervals. The first study concentrated on intervals constructed using empirically determined bandwidths and computed their resulting coverage probabilities. The second study investigated the effects of the iterated smoothed bootstrap on coverage probabilities of both one- and two-sided confidence intervals, with bandwidths fixed in advance. In both studies, we chose $\alpha = 0.05$, $0.1$, $0.9$ and $0.95$, and simulated 1000 random samples of size $n$ from each of three underlying distributions—the standard normal distribution, the double exponential distribution of unit rate and the standard lognormal distribution—in order to estimate the coverage probabilities. The kernels $h, k$ were all taken to be the triangular density function $x \mapsto 1 - |x|$, defined on $[-1, 1]$.

In the first study, we set $q = 0.5$, $n = 15$ and computed the coverage probabilities of $I_{1,\alpha}$ and $I_{3,\alpha}$ constructed using bandwidths determined by the smoothed bootstrap procedure suggested in Section 3. We attempted empirically a wide range of values of $M$ and found that the choice $M = 1.5$ yielded reasonable results under most combinations of $F$ and $\alpha$. We set $M = 1.5$ henceforth. Each interval was constructed using 500 smoothed bootstrap samples, and its coverage probability estimated from 500 outermost bootstrap samples $\mathcal{X}^\circ$. For $I_{1,\alpha}$, the bandwidth $\eta$ for final adoption was searched from 20 evenly spaced values between $0.2n^{-3/8}$ and $2n^{-1/4}$. For $I_{3,\alpha}$, the pilot bandwidths $(\xi, \eta)$ were selected from the $20 \times 20$ grid of evenly spaced values over the rectangle $[0.2n^{-1/2}, 2n^{-3/8}] \times [0.2n^{-1/6}, 2n^{-1/6}]$. Table 2 reports the coverage probabilities of the two intervals, which agree in general with our theoretical result that the Studentized $I_{3,\alpha}$ is more accurate than the non-Studentized $I_{1,\alpha}$. An exception is found for the lognormal data, for which variance estimation is unstable, especially for small samples, and renders $I_{3,\alpha}$ less accurate.

In the second study all four intervals $I_{i,\alpha}$, $i = 1, \ldots, 4$, were constructed for the $q$th population quantiles, for $q = 0.5$ and $0.75$. We also included for reference Beran and Hall's (1993) interpolated interval, denoted by $I_{\mathrm{BH},\alpha}$, and Chen and Hall's (1993) smoothed empirical likelihood intervals, denoted by $I_{\mathrm{EL},\alpha}$ if the interval is not Bartlett-corrected and by $I_{\mathrm{EL}(B),\alpha}$ if it is. Three sample sizes, $n = 15$, 30 and 100, were considered. Each smoothed bootstrap or bootstrap-$t$ interval was constructed using 1000 bootstrap samples. Each iterated interval was constructed using 1500 outer-level and, for each outer-level sample, 1000 inner-level bootstrap samples. For each bootstrap-$t$ interval, we estimated the Studentizing variance by its asymptotic expansion as provided in Sections 2.4 and 2.5, thus avoiding the need for one more level of bootstrapping. Throughout the study, all bandwidths were fixed to be their asymptotically optimal orders for convenience: $\eta = n^{-1/3}$ for $I_{1,\alpha}$; $(\beta, \eta) = (n^{-1/5}, n^{-1/3})$ for $I_{2,\alpha}$; $(\eta, \xi) = (n^{-1/6}, n^{-1/2})$ for $I_{3,\alpha}$; and $(\beta, \eta, \xi) = (n^{-2/19}, n^{-11/57}, n^{-1/2})$ for $I_{4,\alpha}$. For each empirical likelihood interval, we fixed the bandwidth at $n^{-1/2}$, by which Chen and Hall (1993) have produced reasonable results. Tables 3 and 4 summarize the coverage figures for $q = 0.5$ and 0.75, respectively, for the cases $\alpha = 0.05$ and 0.95. Coverage probabilities of two-sided intervals of nominal level 0.9, constructed as $I^2_{\cdot,0.9} = I_{\cdot,0.05} \setminus I_{\cdot,0.95}$, are also reported. Similar findings were obtained for $\alpha = 0.1$ and 0.9, and are therefore omitted from this paper. In the case of two-sided intervals, we estimated also the means and variances of the interval lengths.

We see from Tables 3 and 4 that $I_{2,\alpha}$ is much more accurate than $I_{1,\alpha}$ in most of the cases although the latter is slightly shorter and less variable. This confirms the finite-sample gain acquired by iterating the smoothed bootstrap percentile method. Similar observations are found for the bootstrap-$t$ cases, where the effect of iteration is more notable at the upper end point. However, the degree of improvement of the iterated $I_{4,\alpha}$ over the noniterated $I_{3,\alpha}$ is less remarkable than that achieved by iterating the percentile method, which is not surprising given the very satisfactory coverage already registered by $I_{3,\alpha}$. The coverage figures also demonstrate that $I_{2,\alpha}$ competes closely with $I_{3,\alpha}$ in terms of coverage accuracy. In general $I_{4,\alpha}$ has coverage probabilities comparable to those of the interpolated intervals $I_{\mathrm{BH},\alpha}$. Despite their asymptotically inferior one-sided coverage accuracy, the empirical likelihood intervals $I_{\mathrm{EL},\alpha}$ and $I_{\mathrm{EL}(B),\alpha}$ are found to be slightly more accurate than the bootstrap intervals. Nevertheless the accuracy of the latter improves at a faster rate as $n$ increases compared to that of the empirical likelihood intervals. There is no clear winner in any case, especially for small samples. The two-sided non-Studentized intervals $I^2_{1,0.9}$ and $I^2_{2,0.9}$ are usually shorter and more stable compared to the Studentized intervals $I^2_{3,0.9}$ and $I^2_{4,0.9}$. Note lastly that $I^2_{4,0.9}$ is in general more accurate than $I^2_{\mathrm{EL}(B),0.9}$, although both intervals have coverage errors of orders slightly smaller than $n^{-1}$.

Table 3

*Estimated coverage probabilities of $I_{j,\alpha}$ for $1 - \alpha = 0.05$ ("lower") and $0.95$ ("upper"), and of the $90\%$ two-sided interval $I_{j,0.9}^2$ ("overall"), for $j = 1, 2, 3, 4, \mathrm{BH}, \mathrm{EL}, \mathrm{EL}(B)$*

| Interval | lower | upper | overall | lower | upper | overall | lower | upper | overall |
|---|---|---|---|---|---|---|---|---|---|
| | | $n = 15$ | | | $n = 30$ | | | $n = 100$ | |
| | | | Standard normal data, $N(0,1)$ | | | | | | |
| $I_{1,\alpha}$ | 0.096 | 0.894 | 0.798 | 0.089 | 0.901 | 0.812 | 0.082 | 0.913 | 0.831 |
| $I_{2,\alpha}$ | 0.067 | 0.938 | 0.871 | 0.064 | 0.941 | 0.877 | 0.059 | 0.952 | 0.893 |
| $I_{3,\alpha}$ | 0.057 | 0.932 | 0.875 | 0.054 | 0.935 | 0.881 | 0.058 | 0.935 | 0.877 |
| $I_{4,\alpha}$ | 0.049 | 0.942 | 0.893 | 0.037 | 0.927 | 0.890 | 0.051 | 0.940 | 0.889 |
| $I_{\mathrm{BH},\alpha}$ | 0.046 | 0.952 | 0.906 | 0.051 | 0.950 | 0.899 | 0.046 | 0.950 | 0.904 |
| $I_{\mathrm{EL},\alpha}$ | 0.061 | 0.939 | 0.878 | 0.049 | 0.943 | 0.894 | 0.049 | 0.944 | 0.895 |
| $I_{\mathrm{EL}(B),\alpha}$ | 0.058 | 0.942 | 0.884 | 0.049 | 0.945 | 0.896 | 0.047 | 0.945 | 0.898 |
| | | | Double-exponential data, $\frac{1}{2}\exp(-|x|)$ | | | | | | |
| $I_{1,\alpha}$ | 0.065 | 0.934 | 0.890 | 0.065 | 0.928 | 0.872 | 0.055 | 0.953 | 0.904 |
| $I_{2,\alpha}$ | 0.046 | 0.954 | 0.908 | 0.050 | 0.950 | 0.887 | 0.044 | 0.964 | 0.914 |
| $I_{3,\alpha}$ | 0.044 | 0.955 | 0.908 | 0.049 | 0.951 | 0.899 | 0.040 | 0.952 | 0.916 |
| $I_{4,\alpha}$ | 0.042 | 0.943 | 0.901 | 0.040 | 0.946 | 0.906 | 0.042 | 0.947 | 0.905 |
| $I_{\mathrm{BH},\alpha}$ | 0.031 | 0.957 | 0.926 | 0.059 | 0.951 | 0.892 | 0.046 | 0.957 | 0.911 |
| $I_{\mathrm{EL},\alpha}$ | 0.041 | 0.947 | 0.906 | 0.062 | 0.943 | 0.881 | 0.053 | 0.958 | 0.905 |
| $I_{\mathrm{EL}(B),\alpha}$ | 0.039 | 0.951 | 0.912 | 0.060 | 0.945 | 0.885 | 0.053 | 0.958 | 0.905 |
| | | | Lognormal data, $\exp\{N(0,1)\}$ | | | | | | |
| $I_{1,\alpha}$ | 0.066 | 0.801 | 0.820 | 0.059 | 0.853 | 0.863 | 0.066 | 0.876 | 0.851 |
| $I_{2,\alpha}$ | 0.059 | 0.873 | 0.825 | 0.052 | 0.900 | 0.868 | 0.057 | 0.933 | 0.860 |
| $I_{3,\alpha}$ | 0.046 | 0.888 | 0.865 | 0.045 | 0.918 | 0.881 | 0.055 | 0.927 | 0.871 |
| $I_{4,\alpha}$ | 0.044 | 0.926 | 0.882 | 0.036 | 0.934 | 0.898 | 0.049 | 0.941 | 0.892 |
| $I_{\mathrm{BH},\alpha}$ | 0.047 | 0.956 | 0.909 | 0.051 | 0.950 | 0.899 | 0.046 | 0.951 | 0.905 |
| $I_{\mathrm{EL},\alpha}$ | 0.062 | 0.940 | 0.878 | 0.052 | 0.945 | 0.893 | 0.053 | 0.945 | 0.892 |
| $I_{\mathrm{EL}(B),\alpha}$ | 0.060 | 0.943 | 0.883 | 0.051 | 0.946 | 0.895 | 0.053 | 0.945 | 0.892 |

Sample quantile of interest is $F^{-1}(0.5)$.

**5. Conclusion.** We have examined the asymptotic effects of iterating the smoothed bootstrap on confidence intervals for quantiles, and established the optimal bandwidth orders which minimize the coverage error. Our construction combines two well-known techniques for modifying the conventional bootstrap, smoothing and iteration, in a bootstrap procedure to produce very accurate confidence intervals. Through iterating the smoothed bootstrap, the percentile and bootstrap-$t$ methods yield coverage errors of orders $O(n^{-2/3})$ and $O(n^{-58/57})$, respectively. The latter indeed surpasses all bootstrap methods proposed thus far in the literature as well as Beran and Hall's (1993) interpolated interval and Chen and Hall's (1993) smoothed empirical likelihood intervals with or without Bartlett-correction. The asymptotic gain acquired by iterating the bootstrap in the present context is somewhat

TABLE 4
*Estimated coverage probabilities of $I_{j,\alpha}$ for $1-\alpha = 0.05$ ("lower") and 0.95 ("upper"),*
*and of the 90% two-sided interval $I^2_{j,0.9}$ ("overall"), for $j = 1, 2, 3, 4, \text{BH}, \text{EL}, \text{EL}(B)$*

| Interval | lower | upper | overall | lower | upper | overall | lower | upper | overall |
|---|---|---|---|---|---|---|---|---|---|
| | **$n = 15$** | | | **$n = 30$** | | | **$n = 100$** | | |
| | Standard normal data, $N(0,1)$ | | | | | | | | |
| $I_{1,\alpha}$ | 0.087 | 0.858 | 0.771 | 0.085 | 0.888 | 0.803 | 0.081 | 0.906 | 0.825 |
| $I_{2,\alpha}$ | 0.052 | 0.923 | 0.871 | 0.058 | 0.941 | 0.883 | 0.056 | 0.947 | 0.891 |
| $I_{3,\alpha}$ | 0.044 | 0.885 | 0.841 | 0.061 | 0.926 | 0.865 | 0.062 | 0.949 | 0.887 |
| $I_{4,\alpha}$ | 0.042 | 0.905 | 0.863 | 0.043 | 0.930 | 0.887 | 0.058 | 0.962 | 0.904 |
| $I_{\text{BH},\alpha}$ | 0.039 | 0.979 | 0.940 | 0.053 | 0.950 | 0.897 | 0.046 | 0.944 | 0.898 |
| $I_{\text{EL},\alpha}$ | 0.050 | 0.925 | 0.875 | 0.052 | 0.941 | 0.889 | 0.053 | 0.932 | 0.878 |
| $I_{\text{EL}(B),\alpha}$ | 0.046 | 0.932 | 0.886 | 0.052 | 0.944 | 0.892 | 0.053 | 0.936 | 0.883 |
| | Double-exponential data, $\frac{1}{2}\exp(-|x|)$ | | | | | | | | |
| $I_{1,\alpha}$ | 0.089 | 0.862 | 0.773 | 0.079 | 0.848 | 0.769 | 0.068 | 0.883 | 0.815 |
| $I_{2,\alpha}$ | 0.062 | 0.914 | 0.852 | 0.048 | 0.911 | 0.863 | 0.049 | 0.926 | 0.877 |
| $I_{3,\alpha}$ | 0.044 | 0.890 | 0.846 | 0.065 | 0.927 | 0.862 | 0.060 | 0.938 | 0.878 |
| $I_{4,\alpha}$ | 0.049 | 0.912 | 0.863 | 0.056 | 0.923 | 0.867 | 0.065 | 0.940 | 0.875 |
| $I_{\text{BH},\alpha}$ | 0.035 | 0.974 | 0.939 | 0.057 | 0.961 | 0.904 | 0.049 | 0.944 | 0.895 |
| $I_{\text{EL},\alpha}$ | 0.044 | 0.933 | 0.889 | 0.061 | 0.950 | 0.889 | 0.053 | 0.928 | 0.875 |
| $I_{\text{EL}(B),\alpha}$ | 0.044 | 0.939 | 0.895 | 0.057 | 0.952 | 0.895 | 0.052 | 0.931 | 0.879 |
| | Lognormal data, $\exp\{N(0,1)\}$ | | | | | | | | |
| $I_{1,\alpha}$ | 0.091 | 0.746 | 0.655 | 0.069 | 0.798 | 0.729 | 0.058 | 0.851 | 0.793 |
| $I_{2,\alpha}$ | 0.045 | 0.847 | 0.802 | 0.052 | 0.906 | 0.854 | 0.052 | 0.937 | 0.885 |
| $I_{3,\alpha}$ | 0.049 | 0.844 | 0.795 | 0.042 | 0.899 | 0.857 | 0.053 | 0.928 | 0.875 |
| $I_{4,\alpha}$ | 0.050 | 0.883 | 0.833 | 0.048 | 0.901 | 0.853 | 0.053 | 0.954 | 0.901 |
| $I_{\text{BH},\alpha}$ | 0.039 | 0.979 | 0.940 | 0.053 | 0.950 | 0.897 | 0.046 | 0.945 | 0.899 |
| $I_{\text{EL},\alpha}$ | 0.051 | 0.918 | 0.867 | 0.053 | 0.944 | 0.891 | 0.051 | 0.933 | 0.882 |
| $I_{\text{EL}(B),\alpha}$ | 0.044 | 0.922 | 0.878 | 0.053 | 0.949 | 0.896 | 0.051 | 0.934 | 0.883 |

Sample quantile of interest is $F^{-1}(0.75)$.

nonstandard in the sense that the reduction in coverage error is of an order smaller than $O(n^{-1/2})$ as is commonly the case in regular settings. Table 1 gives a summary of the nonstandard asymptotic improvement effected by smoothed bootstrap iteration.

We have also discussed the effects of using smoothed sample quantiles instead of sample quantiles in the construction of bootstrap intervals or using higher-order instead of second-order kernels for Studentization in the construction of bootstrap-$t$ intervals. Both approaches are shown to yield no asymptotic gain.

Empirical findings of our simulation study agree broadly with the asymptotic theory. Bootstrap-$t$ intervals are in general more accurate than per-centile method intervals; and iterated intervals are more accurate than non-iterated intervals. On the other hand, the percentile method intervals do not require variance estimation for the sample quantile and hence possess the extra advantage of being stable and short compared to the bootstrap-$t$ intervals of the same nominal level.

Bootstrap iteration requires an additional level of bootstrapping, resulting in a computationally more intensive procedure. The apparent computational cost of simulating two batches of outer-level bootstrap samples, the $\mathcal{X}^{\dagger}$'s and the $\mathcal{X}^*$'s, can be alleviated as follows. First simulate one single batch of samples $(\mathcal{Y}^*, \mathcal{W}^*)$'s, where $\mathcal{Y}^* = (Y_1^*, \ldots, Y_n^*)$ and $\mathcal{W}^* = (W_1^*, \ldots, W_n^*)$ denote independent random samples from $F_n$ and $K$, respectively. Combine $\mathcal{Y}^*$ and $\mathcal{W}^*$ using the appropriate bandwidths to form the smoothed bootstrap samples $\mathcal{X}^{\dagger} = (Y_1^* + \eta W_1^*, \ldots, Y_n^* + \eta W_n^*)$ and $\mathcal{X}^* = (Y_1^* + \beta W_1^*, \ldots, Y_n^* + \beta W_n^*)$. Studentization does not pose a computational problem due to the availability of an explicit asymptotic formula for the variance of a sample quantile, which can be readily estimated in practice.

Optimal orders of kernel bandwidths in our construction are specific to the types of intervals being considered. In general, the outer-level smoothed bootstrapping step requires a bandwidth wider than the inner level. For the bootstrap-$t$ method, the bandwidth used for Studentization should be narrower than the bandwidths required by both levels of bootstrapping. The iterated bootstrap typically imposes stricter conditions on our choices of bandwidths. It would be interesting to compare the optimal orders of bandwidths in the quantile problem with those typically recommended for more conventional problems. For instance, density estimation requires a bandwidth of the familiar order $n^{-1/5}$. Under smooth function model settings, asymptotic improvement over the unsmoothed bootstrap can only be effected by a bandwidth of the order $n^{-1/4}$ for the one-sided smoothed bootstrap percentile method interval, and of the order $n^{-1/2}$ for either the two-sided smoothed bootstrap percentile method or bootstrap-$t$ intervals; see Polansky and Schucany (1997).

While we acknowledge the importance of the issue of bandwidth selection and have suggested a bootstrap approach to empirically setting optimal bandwidths, our empirical experience derived from the second simulation study suggests that significant improvement in terms of coverage accuracy can still be acquired, even for small samples, by our simplistic specification of the bandwidths to an arbitrary multiple, which we set as 1 in the above study, of their optimal orders.

## APPENDIX

PROOF OF THEOREM 1.   We follow the proof of Theorem 3.1 in Falk and Janas (1992). Define the distribution function of the standardized sample

$q$th quantile by

$$\Lambda_{n,f}(t) \equiv \mathbb{P}\{n^{1/2}f(F^{-1}(q))[F_n^{-1}(q) - F^{-1}(q)]/\sigma_q \leq t\}, \qquad t \in \mathbb{R}.$$

Under the assumed continuity property of $F$, Reiss [(1989), page 119] established an Edgeworth expansion for $\Lambda_{n,f}$ as

$$
\begin{aligned}
\text{(A.1)} \quad \Lambda_{n,f}(t) = \Phi(t) + n^{-1/2} &\left[ \left( \frac{\delta_n - q}{\sigma_q} + \frac{2(2q-1)}{3\sigma_q} \right) \right. \\
&\left. + \left( \frac{2q-1}{3\sigma_q} + \frac{\sigma_q f'(F^{-1}(q))}{2f^2(F^{-1}(q))} \right) t^2 \right] \phi(t) \\
&+ O(n^{-1}),
\end{aligned}
$$

where $\delta_n = 1 + nq - \lceil nq \rceil$ and $\lceil x \rceil$ denotes the smallest integer greater than or equal to $x$. Close examination of the proof of (A.1) shows that the expansion is actually valid uniformly for $F \in \mathcal{F}_2(\varepsilon, D_1, D_2)$ for any $\varepsilon \in (0, \bar{q})$ and $D_1, D_2 > 0$, a result analogous to Theorem 2.1 of Janas (1993). Such uniform validity carries over to all the expansions which follow, and we presume this fact without mention. We write $\Lambda_n$ for $\Lambda_{n,f}$ if $f$ is the uniform density function over the interval $[0, 1]$. To start with, we require the following bounds for the distances between different versions of sample quantiles, distribution functions and density derivatives:

$$\text{(A.2)} \qquad \hat{F}_{n,\eta}^{-1}(q) - F_n^{-1}(q) = O_p(\eta^2 + n^{-1/2}\eta^{1/2} + n^{-1}),$$

$$\text{(A.3)} \qquad F_n^{-1}(q) - F^{-1}(q) = O_p(n^{-1/2}),$$

$$\text{(A.4)} \qquad \hat{F}_{n,\eta}(q) - F_n(q) = O_p(\eta^2 + n^{-1/2}\eta^{1/2}),$$

$$\text{(A.5)} \qquad \hat{f}_{n,\eta}^{(d)} - f^{(d)} = O_p(\eta^2 + n^{-1/2}\eta^{-d-1/2}) \qquad \text{for } d = 0, 1.$$

For details of the above bounds, see Falk and Janas (1992) for (A.2) and (A.3), Zhou (1997) for (A.4) and Jones (1994) for (A.5).

Now we outline the key steps of the proof. Under the condition that $\eta \propto n^{-\Delta_\eta}$ with $1/4 < \Delta_\eta < 1/2$, straightforward Taylor expansion in conjunction with the bounds (A.2)–(A.5) gives that

$$
\begin{aligned}
\text{(A.6)} \quad &\mathbb{P}\{\hat{G}_{n,\eta}(n^{1/2}[F_n^{-1}(q) - F^{-1}(q)]) < x\} \\
&= \mathbb{P}\left\{ \frac{n^{1/2}}{\sigma_q} \left[ \hat{f}_{n,\eta}(F_n^{-1}(q))(F_n^{-1}(q) - F^{-1}(q)) \right. \right. \\
&\qquad\qquad \left. \left. + \frac{1}{2}\hat{f}'_{n,\eta}(F^{-1}(q))(F_n^{-1}(q) - F^{-1}(q))^2 \right] < \Lambda_n^{-1}(x) \right\} \\
&\quad + O(\eta^2 + n^{-1/2}\eta^{1/2}).
\end{aligned}
$$

Conditioning on the event that $(n^{1/2}/\sigma_q)(F_n^{-1}(q) - F^{-1}(q))f(F^{-1}(q)) = u$, followed by integrating the conditional probability over $u \in \mathbb{R}$ with respect to the distribution of the standardized sample quantile, the probability in (A.6) equals

$$
\text{(A.7)} \quad
\begin{aligned}
\int \mathbb{P}\{ &\hat{f}_{n,\eta}(u_n) u / f(F^{-1}(q)) \\
&+ \tfrac{1}{2} n^{-1/2} \hat{f}'_{n,\eta}(F^{-1}(q)) u^2 \sigma_q / f(F^{-1}(q))^2 \\
&\quad < \Lambda_n^{-1}(x) | F_n^{-1}(q) = u_n \} \Lambda_{n,f}(du),
\end{aligned}
$$

where $u_n = F^{-1}(q) + u n^{-1/2} \sigma_q / f(F^{-1}(q))$. Reiss [[1989], page 119] argues that conditional on $F_n^{-1}(q) = u_n$, $\mathcal{X}$ can be treated as a collection of $n$ independent random variables $Y_1, \ldots, Y_n$, where $\{Y_1, \ldots, Y_{\lceil nq \rceil - 1}\}$ constitutes a random sample from the right truncated density $f(x)/F(u_n)\mathbf{1}\{x < u_n\}$, $Y_{\lceil nq \rceil} = u_n$ and $\{Y_{\lceil nq \rceil + 1}, \ldots, Y_n\}$ is a random sample from the left truncated density $f(x)/(1 - F(u_n))\mathbf{1}\{x > u_n\}$, with $\mathbf{1}\{\cdot\}$ denoting the indicator function. It follows that the conditional probability in (A.7) equals the unconditional probability that $\sum_{i=1}^n T_{n,u,i} < \Lambda_n^{-1}(x)$, for a sum of independent random variables $T_{n,u,i} = T_{n,u,i}(Y_i)$. Let $\mu_{n,u}$ and $\sigma_{n,u}^2$ be the mean and variance of $\sum_{i=1}^n T_{n,u,i}$, respectively. We find by the delta method that $\mu_{n,u} = u\{1 + R_n(u)\}$ and $\sigma_{n,u} = |u| a_n(u)$, where

$$
R_n(u) = n^{-1/2}\left( \frac{2q-1}{2\sigma_q} + \frac{3\sigma_q f'(F^{-1}(q))}{2f^2(F^{-1}(q))} \right) + O(\eta^2 + n^{-2}\eta^{-4})
$$

and

$$
a_n(u) = (n\eta)^{-1/2}\left( \frac{\int k(v)^2 \, dv}{f(F^{-1}(q))} \right)^{1/2} + O(n^{-1/2}\eta^{1/2} + n^{-3/2}\eta^{-5/2}).
$$

Standardizing $\sum_{i=1}^n T_{n,u,i}$ to $S_{n,u} = (\sum_{i=1}^n T_{n,u,i} - \mu_{n,u})/\sigma_{n,u}$ and decomposing the characteristic function of $S_{n,u}$ into factors contributed, respectively, by the partial sums $\sum_{j=1}^{\lceil nq \rceil - 1} T_{n,u,j}$, $T_{n,u,\lceil nq \rceil}$ and $\sum_{j=\lceil nq \rceil + 1}^n T_{n,u,j}$, we can derive an Edgeworth expansion for the distribution of $S_{n,u}$ using standard arguments and rewrite (A.7) as

$$
\text{(A.8)} \quad
\begin{aligned}
\int \mathbb{P}(S_{n,u} &< \vartheta_n(u) + \Theta_n(u)) \Lambda_{n,f}(du) \\
&= \int \Bigg\{ \Phi(\vartheta_n(u)) + \sum_{i=1}^{m_1} \Phi^{(i)}(\vartheta_n(u)) \Theta_n(u)^i / i! \\
&\qquad + \sum_{j=0}^{m_2} (p_{n,u}\phi)^{(j)}(\vartheta_n(u)) \Theta_n(u)^j / j! \Bigg\} \Lambda_{n,f}(du) \\
&\quad + O(\eta^2 + n^{-2}\eta^{-4} + n^{-1/2}\eta^{1/2}) \\
&= \int \{ \Phi(\vartheta_n(u)) + \zeta_n(u) \} \Lambda_{n,f}(du) \\
&\quad + O(\eta^2 + n^{-2}\eta^{-4} + n^{-1/2}\eta^{1/2}) \qquad \text{say,}
\end{aligned}
$$

where $\vartheta_n(u) = (\Lambda_n^{-1}(x) - u)/(|u|a_n(u))$, $\Theta_n(u) = -\text{sign}(u)R_n(u)/a_n(u)$,

$$p_{n,u}(z) = (n\eta)^{-1/2}\left\{\frac{\text{sign}(u)f(F^{-1}(q))\int k(v)^3\,dv}{(f(F^{-1}(q))\int k(v)^2\,dv)^{3/2}}\right\}\left(\frac{z^2-1}{6}\right)\phi(z)$$
$$+ O_p(n^{-3/2}\eta^{-7/2}),$$

and the $m_i$ are chosen such that $(\eta^{-2} + n^2\eta^4 + n^{1/2}\eta^{-1/2})\Theta_n(u)^{m_i+1} \to 0$, for $i = 1, 2$.

Note that $\zeta_n(u)$ is bounded by terms of the form $|\vartheta_n(u)|^j|\Theta_n(u)|^\ell\phi(\vartheta_n(u))$, for positive integers $j, \ell$. For any fixed $u \in (-\infty, 2y] \cup [2y/3, \infty)$ and any fixed $y < 0$, $|\vartheta_n(u)|^j|\Theta_n(u)|^\ell\phi(\vartheta_n(u))$ has order $O(n^{-\lambda})$, for any fixed $\lambda > 0$ and sufficiently large $n$. The same applies to $\Phi(\vartheta_n(u))$ if $u \geq 2y/3$, and to $1 - \Phi(\vartheta_n(u))$ if $u \leq 2y$. Assuming without loss of generality $\Lambda_n^{-1}(x) < 0$ and writing $y = \Lambda_n^{-1}(x)$, it follows from the above bounds and (A.8) that the integral (A.7) equals

$$\begin{aligned}
&\int_y^{2y/3}\Phi(\vartheta_n(u))\Lambda_{n,f}(du) + \int_{2y}^y(\Phi(\vartheta_n(u))-1)\Lambda_{n,f}(du) \\
\text{(A.9)} \quad &+ \int_{2y}^{2y/3}\zeta_n(u)\Lambda_{n,f}(du) + \Lambda_{n,f}(y) + O(\eta^2 + n^{-2}\eta^{-4} + n^{-1/2}\eta^{1/2}) \\
&= I_n + II_n + III_n + \Lambda_{n,f}(y) \\
&\quad + O(\eta^2 + n^{-2}\eta^{-4} + n^{-1/2}\eta^{1/2}) \qquad \text{say.}
\end{aligned}$$

It follows from term-by-term integration with the aid of (A.1) that

$$\begin{aligned}
\text{(A.10)} \quad I_n &= -II_n + O((n\eta)^{-1}) \\
&= (2\pi)^{-1/2}(n\eta)^{-1/2}\left(\frac{\int k(v)^2\,dv}{f(F^{-1}(q))}\right)^{1/2}|y|\Lambda_{n,f}'(y) + o((n\eta)^{-1/2})
\end{aligned}$$

and

$$\begin{aligned}
\text{(A.11)} \quad III_n &= -n^{-1/2}\left(\frac{2q-1}{2\sigma_q} + \frac{3\sigma_q f'(F^{-1}(q))}{2f(F^{-1}(q))^2}\right)y^2\phi(y) \\
&\quad + O(\eta^2 + n^{-2}\eta^{-4} + n^{-1/2}\eta^{1/2}).
\end{aligned}$$

Inverting $\Lambda_n$ and substituting into (A.1), we get

$$\text{(A.12)} \quad \Lambda_{n,f}(y) = x + n^{-1/2}\left(\frac{\sigma_q f'(F^{-1}(q))}{2f(F^{-1}(q))^2}\right)(\Phi^{-1}(x))^2\phi(\Phi^{-1}(x)) + O(n^{-1}).$$

The expansion (1) then follows by combining (A.6), (A.9)–(A.12), putting $x = \alpha$ and taking the complement of the probability. $\square$

PROOF OF THEOREM 2. Arguing as in Janas (1993), we see that the conditions on $k$ and $\beta$ imply that $\mathbb{P}(\hat{F}_{n,\beta} \in \mathcal{F}_2(\varepsilon, D_1, D_2)) = 1 - o(n^{-\lambda})$ for any $\lambda > 0$. It follows that, on substitution of $\hat{F}_{n,\beta}$ for $F$ in (1), the iterated version of the probability has the expansion

$$
\begin{aligned}
(A.13) \quad \hat{J}_{n,\beta,\eta}(x) = {} & x - n^{-1/2}\Big(\frac{2q-1}{2\sigma_q} + \frac{\sigma_q \hat{f}_{n,\beta}(\hat{F}_{n,\beta}^{-1}(q))}{\hat{f}_{n,\beta}(\hat{F}_{n,\beta}^{-1}(q))^2}\Big)\Phi^{-1}(x)^2\phi(\Phi^{-1}(x)) \\
& + O_p(\eta^2 + n^{-2}\eta^{-4} + n^{-1/2}\eta^{1/2}).
\end{aligned}
$$

It follows from the bounds (A.2)–(A.5) that $\hat{J}_{n,\beta,\eta}^{-1}(\alpha)$ differs from the $\alpha$th quantile of the prepivoted root $\hat{G}_{n,\eta}(n^{1/2}[F_n^{-1}(q) - F^{-1}(q)])$ by an order of $O_p(\eta^2 + n^{-2}\eta^{-4} + n^{-1/2}\eta^{1/2} + n^{-1/2}\beta^2 + n^{-1}\beta^{-3/2})$. Theorem 2 then follows by the delta method. $\quad\square$

PROOF OF THEOREM 3. We denote in the sequel by $C_{1,i}(F), C_{2,i}(F), \dots$ generic smooth functions of density derivatives $\{f(F^{-1}(q)), f'(F^{-1}(q)), \dots, f^{(i)}(F^{-1}(q))\}$, for each $i = 0, 1, \dots$. In cases where $C_{j,i}(F)$ assumes the form of a polynomial in a variable $x$, we write $C_{j,i}(F) = C_{j,i}(F)(x)$. As in (A.7), we define $u_n = F^{-1}(q) + un^{-1/2}\sigma_q/f(F^{-1}(q))$ and write

$$
(A.14) \quad K_{n,\xi}(x) = \int \mathbb{P}\{\hat{f}_{n,\xi}(u_n)u/f(F^{-1}(q)) < x | F_n^{-1}(q) = u_n\}\Lambda_{n,f}(du).
$$

Arguing as in the proof of Theorem 1, we show that the conditional distribution of the standardized form of $\hat{f}_{n,\xi}(u_n)u/f(F^{-1}(q))$, given $F_n^{-1}(q) = u_n$, has an Edgeworth expansion

$$
\begin{aligned}
(A.15) \quad & \Phi(y) + (n\xi)^{-1/2}C_{1,0}(F)(y) \\
& + (n\xi)^{-1}C_{2,0}(F)(y) + n^{-1/2}\xi^{1/2}C_{3,0}(F)(y) \\
& + (n\xi)^{-3/2}C_{4,0}(F)(y) + n^{-1}\xi^{-1/2}C_{5,1}(F)(y) \\
& + O_p(n^{-1} + n^{-3/2}\xi^{-1} + (n\xi)^{-2} + n^{-1/2}\xi^{3/2}).
\end{aligned}
$$

Plugging (A.15) into (A.14) and splitting the integral as in (A.9), we see in the present context that

$$
\begin{aligned}
(A.16) \quad I_n + II_n = {} & (n\xi)^{-1}C_{6,0}(F)(x) \\
& + (n\xi)^{-3/2}C_{7,0}(F)(x) + (n\xi)^{-2}C_{8,0}(F)(x) \\
& + n^{-1}C_{9,0}(F)(x) + n^{-3/2}\xi^{-1}C_{10,1}(F)(x) \\
& + O(n^{-1}\xi^{1/2} + n^{-2}\xi^{-3/2} + (n\xi)^{-5/2})
\end{aligned}
$$

and

$$
\begin{aligned}
(A.17) \quad III_n = {} & (n\xi)^{-1}C_{11,0}(F)(x) + n^{-1/2}C_{12,1}(F)(x) \\
& + (n\xi)^{-3/2}C_{13,0}(F)(x) + n^{-1}\xi^{-1/2}C_{14,1}(F)(x) \\
& + n^{-1/2}\xi^{1/2}C_{15,1}(F)(x) + n^{-1}C_{16,2}(F)(x) \\
& + n^{-3/2}\xi^{-1}C_{17,1}(F)(x) + \xi^2 C_{18,2}(F)(x) \\
& + n^{-1/2}\xi C_{19,1}(F)(x) + (n\xi)^{-2}C_{20,0}(F)(x) \\
& + O(n^{-3/2}\xi^{-1/2} + (n\xi)^{-5/2} + n^{-1/2}\xi^{3/2} \\
& \quad + n^{-2}\xi^{-3/2} + n^{1/2}\xi^{9/2} + \xi^{5/2} + n^{3/2}\xi^{19/2}).
\end{aligned}
$$

It follows by noting (A.1) and combining (A.16) and (A.17) that

$$
\begin{aligned}
K_{n,\xi}(x) ={} & \Phi(x) + (n\xi)^{-1}C_{21,0}(F)(x) \\
& + n^{-1/2}C_{22,1}(F)(x) + (n\xi)^{-3/2}C_{23,0}(F)(x) \\
& + n^{-1}\xi^{-1/2}C_{24,1}(F)(x) \\
& + n^{-1/2}\xi^{1/2}C_{25,1}(F)(x) + n^{-3/2}\xi^{-1}C_{26,1}(F)(x) \\
& + (n\xi)^{-2}C_{27,0}(F)(x) + n^{-1}C_{28,2}(F)(x) + \xi^2 C_{29,2}(F)(x) \\
& + n^{-1/2}\xi C_{30,1}(F)(x) \\
& + O(n^{-2}\xi^{-3/2} + (n\xi)^{-5/2} + n^{1/2}\xi^{9/2} + \xi^{5/2} + n^{3/2}\xi^{19/2}).
\end{aligned}
\tag{A.18}
$$

Note that (A.18) extends the expansion given in Theorem 2.1 of Janas (1993) by including higher-order terms. Arguments similar to those provided by Janas (1993) can be used to show that (A.18) holds uniformly in $F \in \mathcal{F}_3(\varepsilon, D_1, D_2)$. The assumptions (A1)–(A3) and that $\Delta_\eta \in (0, 1/5)$ imply that $\hat{F}_{n,\eta} \notin \mathcal{F}_3(\varepsilon, D_1, D_2)$ with negligible probability. An asymptotic expansion similar to (A.18) thus holds for $\hat{K}_{n,\eta,\xi}(x)$, with $F$ substituted by $\hat{F}_{n,\eta}$. Standard Taylor expansion together with bounds (A.2)–(A.5), with (A.5) strengthened to include cases $d = 2, 3$, yields an expansion for the difference

$$
\begin{aligned}
\hat{K}_{n,\eta,\xi}^{-1}(x) &- K_{n,\xi}^{-1}(x) \\
={} & \hat{\kappa}_{n,\eta,\xi}(x) + O_p((n\xi)^{-1}\eta^4 + n^{-1/2}\eta^4 + n^{-3/2}\xi^{-1}\eta^{1/2} + n^{-1}\eta^{1/2} \\
& \qquad + (n\xi)^{-5/2} + \xi^{5/2} + n\xi^5 + n^{-3/2}\eta^{-5/2} \\
& \qquad + n^{-1/2}\xi^2\eta^{-5/2} + n^{-1/2}\xi^{1/2}\eta^2 + (n\xi)^{-3/2}\eta^2),
\end{aligned}
$$

where

$$
\begin{aligned}
\hat{\kappa}_{n,\eta,\xi}(x) ={} & \{(n\xi)^{-1}C_{31,0}(F)(x) + n^{-1/2}C_{32,1}(F)(x) \\
& \quad + (n\xi)^{-3/2}C_{33,0}(F)(x) + n^{-1}\xi^{-1/2}C_{34,1}(F)(x) \\
& \qquad\qquad + n^{-1/2}\xi^{1/2}C_{35,1}(F)(x)\} \\
& \times [\hat{f}_{n,\eta}(F_n^{-1}(q)) - f(F^{-1}(q))] \\
& + \{n^{-1/2}C_{36,1}(F)(x) + n^{-1}\xi^{-1/2}C_{37,1}(F)(x) \\
& \qquad\qquad + n^{-1/2}\xi^{1/2}C_{38,1}(F)(x)\} \\
& \times [\hat{f}'_{n,\eta}(F_n^{-1}(q)) - f'(F^{-1}(q))] \\
& + \{(n\xi)^{-1}\eta^2 C_{39,1}(F)(x) + n^{-1/2}\eta^2 C_{40,1}(F)(x)\}\hat{f}'_{n,\eta}(F_n^{-1}(q)) \\
& + n^{-1/2}\eta^2 C_{41,1}(F)(x)\hat{f}''_{n,\eta}(F_n^{-1}(q)).
\end{aligned}
$$

Note, by conditioning and integrating as in (A.14), that

$$
\begin{aligned}
(A.19) \quad & \mathbb{P}\{\hat{K}_{n,\eta,\xi}(n^{1/2}[(F_n^{-1}(q) - F^{-1}(q))/\hat{s}_\xi]) \geq \alpha\} \\
& = \int \mathbb{P}\{\hat{f}_{n,\xi}(u_n)u/f(F^{-1}(q)) \\
& \quad - \hat{\kappa}_{n,\eta,\xi}(\alpha) \geq K_{n,\xi}^{-1}(\alpha)|F_n^{-1}(q) = u_n\}\Lambda_{n,f}(du) \\
& \quad + O((n\xi)^{-1}\eta^4 + n^{-1/2}\eta^4 + n^{-3/2}\xi^{-1}\eta^{1/2} \\
& \quad + n^{-1}\eta^{1/2} + (n\xi)^{-5/2} + \xi^{5/2} + n\xi^5 + n^{-3/2}\eta^{-5/2} \\
& \quad + n^{-1/2}\xi^2\eta^{-5/2} + n^{-1/2}\xi^{1/2}\eta^2 + (n\xi)^{-3/2}\eta^2).
\end{aligned}
$$

Consider first the integral

$$
\tilde{K}_{n,\eta,\xi}(y) = \int \mathbb{P}\{\hat{f}_{n,\xi}(u_n)u/f(F^{-1}(q)) - \hat{\kappa}_{n,\eta,\xi}(\alpha) \leq y|F_n^{-1}(q) = u_n\}\Lambda_{n,f}(du).
$$

Proceeding, with lengthy algebra, as in establishing the Edgeworth expansion for the conditional probability in (A.7), we see that the conditional distributions of the standardized $\hat{f}_{n,\xi}(u_n)u/f(F^{-1}(q))$ and $\hat{f}_{n,\xi}(u_n)u/f(F^{-1}(q)) - \hat{\kappa}_{n,\eta,\xi}(\alpha)$, given that $F_n^{-1}(q) = u_n$, have the same Edgeworth expansion up to $O_p((n\xi)^{-2} + n^{-1/2}\xi^{3/2} + n^{-3/2}\xi^{-1/2}\eta^{-1} + n^{-1}\xi^{1/2}\eta^{-1} + n^{-1}\xi^{3/2}\eta^{-3})$. It follows, by lengthy algebra again, that the extra term $\hat{\kappa}_{n,\eta,\xi}(\alpha)$ in (A.19) contributes only to the $III_n$ component of the integral (A.14) through $\mathbb{E}\hat{\kappa}_{n,\eta,\xi}(\alpha)$, up to order

$$
\begin{aligned}
& O((n\xi)^{-5/2} + n\xi^5 + \xi^{5/2} + n^{-1/2}\eta^4 \\
& \quad + n^{-1}\eta^{1/2} + n^{-3/2}\eta^{-5/2} + n^{-3/2}\xi^{-1/2}\eta^{-1} \\
& \quad + (n\xi)^{-1}\eta^4 + n^{-1}\xi^{-1/2}\eta^2 + (n\eta)^{-1}\xi^{1/2} + n^{-1/2}\xi^{1/2}\eta^2 \\
& \quad + n^{-3/2}\xi^{-1}\eta^{1/2} + n^{-1/2}\xi^2\eta^{-5/2} + (n\xi)^{-3/2}\eta^2).
\end{aligned}
$$

Explicit expansion of the contribution of $\hat{\kappa}_{n,\eta,\xi}(\alpha)$ yields that $\tilde{K}_{n,\xi}(y)$ differs from $K_{n,\xi}(y)$ by

$$
\begin{aligned}
& (n\eta)^{-1}C_{42,1}(F)(y) + (n\xi)^{-1}\eta^2 C_{43,2}(F)(y) + n^{-1/2}\eta^2 C_{44,3}(F)(y) \\
& \quad + n^{-1}C_{45,2}(F)(y) + n^{-3/2}\xi^{-1}C_{46,1}(F)(y)
\end{aligned}
$$

up to the above order, uniformly in $F \in \mathcal{F}_4(\varepsilon, D_1, D_2)$. Theorem 3 now follows by putting $y = K_{n,\xi}^{-1}(\alpha)$ and taking the complement.  □

PROOF OF THEOREM 4. As in the proof of Theorem 2, when $0 < \Delta_\beta < 1/7$, we have $\hat{F}_{n,\beta} \notin \mathcal{F}_4(\varepsilon, D_1, D_2)$ with negligible probability. Application of the bounds (A.2)–(A.4) and extending (A.5) to cases $d = 2, 3$ establish that the bootstrap quantile $\hat{L}_{n,\beta,\eta,\xi}^{-1}(\alpha)$ and the true quantile of $\hat{K}_{n,\eta,\xi}(n^{1/2}[(F_n^{-1}(q) - F^{-1}(q))/\hat{s}_\xi])$ differ in probability by the error term as specified in (3). Theorem 4 then follows by the delta method.  □

PROOF OF THEOREM 5. Write $\Psi_{n,\zeta}(t) \equiv G_{n,\zeta}(t\sigma_q/f(F^{-1}(q)))$ and $\hat{\Psi}_{n,\eta,\zeta}(t) \equiv \hat{G}_{n,\eta,\zeta}(t\sigma_q/\hat{f}_{n,\eta}(\hat{F}_{n,\eta}^{-1}(q)))$ for $t \in \mathbb{R}$. Standard Edgeworth expansion gives that

$$
\begin{aligned}
(A.20) \quad \Psi_{n,\zeta}(t) = {} & \Phi(t) + \{n^{-1/2}E_{1,1}(F) + \zeta E_{2,0}(F) + n^{1/2}\zeta^2 E_{3,1}(F)\}\phi(t) \\
& + O(n^{-1} + n^{3/2}\zeta^4)
\end{aligned}
$$

uniformly in $F \in \mathcal{F}_2(\varepsilon, D_1, D_2)$ for any $\varepsilon \in (0, \bar{q})$ and $D_1, D_2 > 0$, where $E_{j,i}(F)$ denotes a smooth function of the density derivatives $\{f(F^{-1}(q)), \ldots, f^{(i)}(F^{-1}(q))\}$ for each $j = 1, 2, 3$. On substitution of $\hat{F}_{n,\eta}$ for $F$ in (A.20) and using (A.2)–(A.5), we have

$$
\begin{aligned}
\hat{\Psi}_{n,\eta,\zeta}^{-1}(\alpha) = {} & \Psi_{n,\zeta}^{-1}(\alpha) \\
& + O_p(n^{-1} + n^{3/2}\zeta^4 + \eta^2\zeta + n^{-1/2}\eta^{-1/2}\zeta + n^{1/2}\eta^2\zeta^2 + \eta^{-3/2}\zeta^2).
\end{aligned}
$$

We obtain by the delta method that the coverage probability of $I_{1,\alpha}^\kappa$ equals

$$
\begin{aligned}
(A.21) \quad & \mathbb{P}\{n^{1/2}(\tilde{F}_{n,\zeta}^{-1}(q) - F^{-1}(q)) - \Psi_{n,\zeta}^{-1}(\alpha)\Delta \geq \Psi_{n,\zeta}^{-1}(\alpha)\sigma_q/f(F^{-1}(q))\} \\
& + O(n^{-1} + n^{3/2}\zeta^4 + \eta^2\zeta + n^{-1/2}\eta^{-1/2}\zeta + n^{1/2}\eta^2\zeta^2 + \eta^{-3/2}\zeta^2),
\end{aligned}
$$

where $\Delta = \sigma_q[1/\hat{f}_{n,\eta}(\hat{F}_{n,\eta}^{-1}(q)) - 1/f(F^{-1}(q))]$. Replacing $\hat{F}_{n,\eta}$ by $\tilde{F}_{n,\zeta}$ in (A.2), using (A.3), (A.5) and noting the bounds assumed on $\eta, \zeta$, we have that the joint cumulants of $(n^{1/2}(\tilde{F}_{n,\zeta}^{-1}(q) - F^{-1}(q)), \Delta)$ differ from those of $(n^{1/2}(F_n^{-1}(q) - F^{-1}(q)), \Delta)$ by at most $O(n^{-25/48})$. An expansion analogous to (A.21) holds for the interval $I_{1,\alpha}$ with the definition of $\Delta$ unchanged. Comparison with (1) then implies that the joint cumulants of $(n^{1/2}(F_n^{-1}(q) - F^{-1}(q)), \Delta)$ contribute a term of precise order $n^{-1/2}$ to the coverage error of $I_{1,\alpha}$. It follows that the Edgeworth expansions for the distributions of $n^{1/2}(\tilde{F}_{n,\zeta}^{-1}(q) - F^{-1}(q))$ and $n^{1/2}(\tilde{F}_{n,\zeta}^{-1}(q) - F^{-1}(q)) - \Psi_{n,\zeta}^{-1}(\alpha)\Delta$ also differ by an order of $n^{-1/2}$ precisely, so that, by (A.21), the coverage probability of $I_{1,\alpha}^\kappa$ has the expansion stated in (4). $\square$

DEPARTMENT OF STATISTICS AND
ACTUARIAL SCIENCE
THE UNIVERSITY OF HONG KONG
POKFULAM ROAD
HONG KONG
E-MAIL: hohs@graduate.hku.hk
E-MAIL: smslee@hkusua.hku.hk